# OPTIMAL THIRD ROOT ASYMPTOTIC BOUNDS IN THE STATISTICAL ESTIMATION OF THRESHOLDS

By Franz Merkl and Leila Mohammadi

*University of Munich and University of Leiden*

This paper is concerned with estimating the intersection point of two densities, given a sample of both of the densities. This problem arises in classification theory. The main results provide lower bounds for the probability of the estimation errors to be large on a scale determined by the inverse cube root of the sample size. As corollaries, we obtain probabilistic bounds for the prediction error in a classification problem. The key to the proof is an entropy estimate. The lower bounds are based on bounds for general estimators, which are applicable in other contexts as well. Furthermore, we introduce a class of optimal estimators whose errors asymptotically meet the border permitted by the lower bounds.

## 1. Introduction.

1.1. *Motivation and origin of the problem.* In this paper we derive lower bounds for the probability of large errors of some estimators to occur. Let $\mathcal{P}$ be a class of probability measures on a measurable space $(\Omega, \mathcal{A})$, and $a : \mathcal{P} \to \mathbb{R}$ be a parameter. Consider an i.i.d. random sample $Z_1, \ldots, Z_n$ from $P \in \mathcal{P}$ and an estimator $\hat{a}_n(Z_1, \ldots, Z_n)$, $\hat{a}_n : \Omega^n \to \mathbb{R}$, for $a$. We are interested in the asymptotic behavior of $\hat{a}_n$ as $n \to \infty$.

In the theory of empirical processes one usually considers a deterministic loss function whose minimizer over a particular class is equal to or close to the parameter. Under some technical assumptions, if the loss is differentiable with respect to the parameter, the empirical risk minimizers converge to the parameter with the rate $|\hat{a}_n - a(P)| = O_P(n^{-1/2})$ as the sample size $n$ grows to $\infty$; see van de Geer [19] and van der Vaart and Wellner [20].

Kim and Pollard [9] establish a new functional central limit theorem for empirical processes. They describe an interesting class of asymptotic problems where the estimators converge at a rate different from $O_P(n^{-1/2})$ to limit distributions.

An important noncontinuous loss function, frequently used in the theory of classification, is the indicator loss function. Let us first describe a general view of classification (or learning theory). We formulate the simplest case, which is a two-class problem. Assume that we have two distributions on a result space $\mathcal{X}$, labeled by $Y=1$ and $Y=-1$. The values of $Y$ are called "labels" or "natures." Take an observation $X$ from a mixture of the two distributions. It is sometimes called a "feature." The problem is to predict the unknown nature $Y$ of a feature $X$. Suppose we have $n$ i.i.d. copies $(X_i,Y_i)$, $i=1,\dots,n$, of a realization $(X,Y)$, having an unknown probability distribution $P$. A classifier $h$ is a measurable function $h:\mathcal{X}\to\{\pm1\}$. (Here, we do not consider more general $[-1,1]$-valued classifiers.) The realization $(X,Y)$ is called misclassified by the classifier $h$ if $h(X)\neq Y$. We take the deterministic loss function $(x,y)\mapsto \mathbb{1}\{h(x)\neq y\}$. For $\mathcal{X}\subseteq\mathbb{R}$ and features $X$ with a continuous distribution (at least close to a point), Mohammadi and van de Geer [13] apply this setup to the case where the classifier $h$ is varied over the class $\mathcal{H}=(h_a)_{a\in\mathbb{R}}$, where $h_a(x):=1$ for $x\geq a$ and $h_a(x):=-1$ for $x<a$. Let

$$f_P(x,y)=f_P^+(x)\mathbb{1}\{y=1\}+f_P^-(x)\mathbb{1}\{y=-1\},\qquad (x,y)\in\Omega=\mathcal{X}\times\{\pm1\}, \tag{1.1}$$

denote the joint density of $(X,Y)$, that is, the density of $P$ with respect to some reference measure $\lambda\otimes$ (counting measure).

Let us assume that there is a unique point $a(P)$ at which $f_P^+-f_P^-$ changes its sign, and that this sign change is from "$-$" to "$+$." Then, to minimize the risk $P[h(X)\neq Y]$, it suffices to restrict the classifier $h$ to the class $\mathcal{H}$, since one has in this case

$$\inf_{\text{classifiers } h} P[h(X)\neq Y]=\inf_{a\in\mathbb{R}} P[h_a(X)\neq Y]=P[h_{a(P)}(X)\neq Y]. \tag{1.2}$$

The *Bayes rule* in this case corresponds to the threshold $a(P)=\arg\min_{a\in\mathbb{R}} L_P(a)$, where $L_P(a):=P(h_a(X)\neq Y)$ denotes the prediction error. A natural choice for an estimator of $a$ is the threshold $\hat{a}_n=\arg\min_{a\in\mathbb{R}}\hat{P}_n[h_a(X)\neq Y]$ that minimizes the classification error in the sample, where $\hat{P}_n:=\sum_{i=1}^n \delta_{(X_i,Y_i)}/n$ denotes the empirical distribution of the sample. Strictly speaking, here the "arg min" is not unique, but one may take any (measurable) choice. In the theory of classification this is called empirical risk minimization. Mohammadi and van de Geer [13] invoke the theory of Kim and Pollard [9] to get the rate $O_P(n^{-1/3})$ of $\hat{a}_n$ under some conditions. Under monotonicity assumptions, it is shown that $\hat{a}_n$ is a nonparametric maximum likelihood estimator, and that $n^{-1/3}(\hat{a}_n - a(P))$ converges in distribution to a continuous random variable. For more background information about empirical risk minimization in classification, see, for instance, [7], [16] and [11].

1.2. *Statement of the problem and results.* In this paper we address the following question: Is there any sequence of estimators $(\hat{a}_n : \Omega^n \to \mathbb{R})_{n\in\mathbb{N}}$ which converges to $a(P)$ with a rate faster than $O_P(n^{-1/3})$? Under some assumptions, to be specified below, the answer is *no*.

Let us introduce the class $\mathcal{P}$ of probability measures $P$ that we consider. We assume that the feature $X$ takes values in the unit interval $\mathcal{X} = [0,1]$.

Let $\tilde{\mathcal{P}}$ denote the set of all probability distributions $P = f_P[\lambda_{[0,1]} \otimes (\text{counting measure})]$ on $\Omega = [0,1] \times \{\pm 1\}$ (with the Borel $\sigma$-field) with $f_P \in C^1([0,1] \times \{\pm 1\})$. Here, $\lambda_{[0,1]}$ denotes the Lebesgue measure on $[0,1]$. [This particular choice of the reference measure—at least locally— plays a role in some technical estimates, e.g., in the basic entropy bound (4.13) below.]

We endow $\tilde{\mathcal{P}}$ with the metric $d$, given by

$$d(P,Q) := \|f_P - f_Q\|_\infty + \|\partial_1 f_P - \partial_1 f_Q\|_\infty, \tag{1.3}$$

where $\partial_1$ denotes the derivative with respect to the first argument.

Let $\mathcal{P} \subseteq \tilde{\mathcal{P}}$ denote the set of all $P \in \tilde{\mathcal{P}}$, such that $f_P^+ := f_P(\cdot, 1)$ and $f_P^- := f_P(\cdot, -1)$ have a unique intersection point $a(P)$, and this intersection point is contained in the open interval $(0,1)$, and the intersection is transversal with a specified orientation,

$$f_P^+(a(P)) = f_P^-(a(P)), \qquad (f_P^+)'(a(P)) > (f_P^-)'(a(P)). \tag{1.4}$$

We endow $\mathcal{P}$ with the topology induced by the metric $d$.

For our results it is essential to have at least some control on the derivative of $f_P$, which is reflected by the choice (1.3) of the metric $d$.

Here, we present the main results of this paper. The first theorem considers the estimation error on the *critical scale* const $\cdot n^{-1/3}$, uniformly over (small) open subsets of $\mathcal{P}$.

THEOREM 1.1. *Let $\mathcal{U} \subseteq \mathcal{P}$ be a nonempty open set. Then there is $c_1 = c_1(\mathcal{U}) > 0$, such that for every $\delta \in (0, 1/4]$ and for every sequence of estimators $\hat{a}_n : \Omega^n \to \mathbb{R}$, $n \in \mathbb{N}$ one has*

$$\liminf_{n\to\infty} \sup_{P\in\mathcal{U}} P^n[n^{1/3}|\hat{a}_n - a(P)| > T] \geq \delta, \tag{1.5}$$

*where $T = T(\mathcal{U}, \delta) := c_1 |\log(11\delta)|^{1/3}$.*

Unlike Theorem 1.1, the following theorem considers the asymptotics of the estimation error *point-wise*, that is, it takes the limit as $n \to \infty$ before taking a supremum over open sets $\mathcal{U} \subseteq \mathcal{P}$. To motivate this order of taking limits, consider the following game of a statistician against "nature." Nature chooses just one $P \in \mathcal{P}$, unknown to the statistician. The statistician may choose various sample sizes $n$, and she or he obtains a certain rate of convergence of the estimators as $n \to \infty$ for this fixed, given $P$. Thus, examining the limit $n \to \infty$ for fixed, but arbitrary $P$ may contain at least as relevant information as taking $\liminf_{n\to\infty} \sup_{P\in\mathcal{U}}$. The following theorem does not examine the critical scale $n^{-1/3}$; it rather works with a smaller scale $1/\beta_n \ll n^{-1/3}$. The reason for this is explained below.

THEOREM 1.2. *Let $(\beta_n)_{n\in\mathbb{N}}$ be a sequence of positive numbers with $\lim_{n\to\infty} n^{-1/3}\beta_n = \infty$. Then, for all nonempty open sets $\mathcal{U} \subseteq \mathcal{P}$ and for all sequences of estimators $\hat{a}_n : \Omega^n \to \mathbb{R}$, $n \in \mathbb{N}$, one has*

$$\sup_{P\in\mathcal{U}} \limsup_{n\to\infty} P^n[\beta_n|\hat{a}_n - a(P)| > 1] \geq 1/4. \tag{1.6}$$

Theorem 1.2 states that, independently of how small our statistical model $\mathcal{U}$ is, we always find a particular model $P$ in this class such that the estimation error for this particular model will be with positive probability asymptotically larger than any given scale smaller than $n^{-1/3}$. The proof of this theorem uses Baire's theorem. A related argument, concerning the equicontinuity and the consistency of substitution estimators with values in a metric space, is presented in [15].

In contrast to Theorem 1.1, Theorem 1.2 does not consider the critical scale const $\cdot n^{-1/3}$. Indeed, its claim (1.6) breaks down on this critical scale. This is the content of the following theorem.

THEOREM 1.3. *There is a family of estimators $(\hat{a}_{n,L} : \Omega^n \to \mathbb{R})_{n\in\mathbb{N},L>0}$ with the following property: For all $P \in \mathcal{P}$, there is a neighborhood $\mathcal{N} \subseteq \mathcal{P}$ of $P$, such that for all $T > 0$ one has*

$$\lim_{L\to\infty} \sup_{Q\in\mathcal{N}} \limsup_{n\to\infty} Q^n[n^{1/3}|\hat{a}_{n,L} - a(Q)| > T] = 0. \tag{1.7}$$

Such estimators $\hat{a}_{n,L}$ are explicitly described in Section 5 below. Speaking very roughly, one estimates the density $f_P$ in a certain neighborhood of size $Ln^{-1/3}$ of a first approximation of $a(Q)$ using regression lines.

The following corollaries translate the asymptotic bounds in Theorems 1.1, 1.2 and 1.3 to bounds for the rate of convergence of the prediction error $L_P(\hat{a}_n)$ to the optimal value $L_P(a(P))$.

COROLLARY 1.4. *Under the conditions of Theorem 1.1, one has*

$$\liminf_{n\to\infty} \sup_{P\in\mathcal{U}} P^n[L_P(\hat{a}_n) > L_P(a(P)) + Sn^{-2/3}] \geq \delta, \tag{1.8}$$

*where $S = S(\mathcal{U}, \delta) := c_2|\log(11\delta)|^{2/3}$ with a constant $c_2 = c_2(\mathcal{U}) > 0$.*

Corollary 1.5. *Under the conditions of Theorem 1.2, in particular, for* $\beta_n^{-2} = o(n^{-2/3})$, *one has*

$$\sup_{P \in \mathcal{U}} \limsup_{n \to \infty} P^n[L_P(\hat{a}_n) > L_P(a(P)) + c_3 \beta_n^{-2}] \geq 1/4, \tag{1.9}$$

*with a constant* $c_3 = c_3(\mathcal{U}) > 0$.

Corollary 1.6. *Take a family of estimators* $(\hat{a}_{n,L})_{n \in \mathbb{N}, L > 0}$ *that fulfills the claim of Theorem 1.3. Then, for all* $P \in \mathcal{P}$, *there is a neighborhood* $\mathcal{N} \subseteq \mathcal{P}$ *of* $P$, *such that, for all* $T > 0$, *one has*

$$\lim_{L \to \infty} \sup_{Q \in \mathcal{N}} \limsup_{n \to \infty} Q^n[L_Q(\hat{a}_n) > L_Q(a(Q)) + T n^{-2/3}] = 0. \tag{1.10}$$

1.3. *Discussion and comparison to other results.* Let us discuss our results and compare them with some previous results on lower bounds in classification and regression.

The paper [10] by Mammen and Tsybakov views the classification problem as the estimation problem of a set $V$. The authors consider the case that the region $V$ has a smooth boundary or belongs to another nonparametric class of sets. They show that the empirical risk minimizers achieve the optimal rates for estimation of $V$ and optimal rates of convergence for Bayes risks.

It is interesting to compare our Theorem 1.1 with Theorem 3 in Mammen and Tsybakov's paper [10], in particular, with formula (22) there. The setup in the paper [10] is much more general than ours. It differs from the one in Theorem 1.1, even when one specializes it to our one-dimensional setup and to the special classifiers $h_a$. More specifically, this specialization yields, for all $p \geq 1$,

$$\liminf_{n \to \infty} \sup_{P \in \mathcal{F}_{\text{frag}}} n^{p/3} E_{P^n}[|\hat{a}_n - a(P)|^p] > 0, \tag{1.11}$$

instead of our claim (1.5), where the class of distributions $\mathcal{F}_{\text{frag}}$ specified in the reference is not as small as our open set $\mathcal{U}$.

Let us compare the estimators $\hat{a}_{n,L}$ in Theorem 1.3 and Corollary 1.6 with the empirical risk minimizers $\hat{a}_n$, which are examined in the paper [13] by Mohammadi and van de Geer. For the empirical risk minimizers $\hat{a}_n$, one has

$$\sup_{Q \in \mathcal{N}} \limsup_{n \to \infty} Q^n[n^{1/3}|\hat{a}_n - a(Q)| > T] > 0 \tag{1.12}$$

for all $T < \infty$; more details are given in Theorem 2.2 in [13]. The empirical risk minimizers $\hat{a}_n$ may be well applicable for larger classes of distributions than $\mathcal{P}$, where our results may not apply. Our intention behind Theorem 1.3 is mainly to show that Theorem 1.2 is optimal. However, comparing (1.12) with (1.7), one sees that, for large $L$, $\hat{a}_{n,L}$ is an improvement over $\hat{a}_n$, at least asymptotically in the limit as $n \to \infty$. Thus, from a practical point of view, we suggest use of the estimators $\hat{a}_{n,L}$ instead of $\hat{a}_n$ whenever one suspects the regularity conditions imposed in our model are reasonable in an application at hand.

Roughly speaking, the improvement is obtained by using information about the empirical distribution in the neighborhood of the estimator $\hat{a}_n$. The scaling parameter $L$ is used to determine the size of this neighborhood. More specifically, one estimates the unknown densities close to $\hat{a}_n$ using regression lines. For a given sample size $n$, it might not make sense to take $L$ too large to get a good estimator $\hat{a}_{n,L}$, due to the order of limits "$\lim_{L\to\infty} \ldots \limsup_{n\to\infty}$" in (1.7).

Donoho and Liu [6] consider estimating a functional $T(F)$ of an unknown distribution $F \in \mathcal{F}$ with some class of distributions $\mathcal{F}$. They compute the modulus of continuity $\omega(\varepsilon)$ of $T$ with respect to Hellinger distance in certain cases. For a well-behaved loss function $l(t)$, they show that if $T$ is linear and $\mathcal{F}$ is convex, then $\inf_{T_n} \sup_{F\in\mathcal{F}} E_F(l(T_n - T(F)))$ is equivalent to $l(\omega(n^{-1/2}))$ within constants. The same conclusion is drawn for three cases of nonlinear functionals: estimating the rate of decay of a density, estimating the mode and robust nonparametric regression. Our case, estimating the intersection point of two densities, is a different case. However, it gets the modulus of continuity $\omega(\varepsilon) = \varepsilon^{2/3}$ for $l(t) := |t|$ and therefore, $\omega(n^{-1/2}) = n^{-1/3}$, which coincides with the optimal rate.

The general estimates for lower bounds, presented in Section 3 below, can also be applied to higher-dimensional problems. This will be shown in a forthcoming paper.

Let us briefly review some further known results which are vaguely related to the facts proven in this paper.

Let $(X,Y), (X_1,Y_1), (X_2,Y_2), \ldots$ be independent identically distributed $\mathbb{R}^d \times \mathbb{R}$ random variables with $E(Y^2) < \infty$. In a regression problem, Stone [17] showed that for a class of distributions and for a class of regression functions which are $p$ times continuously differentiable, the optimal lower rate of convergence is $n^{-2p/(2p+d)}$.

Antos, Györfi and Kohler [2] showed that there exist individual lower bounds on the rate of convergence of nonparametric regression estimates which are arbitrarily close to Stone's minimax lower bounds.

In classification Antos and Lugosi [3] showed that for several natural concept classes (classes of subsets of $\mathcal{X}$, the domain of $X$), including the class of linear half-spaces, there exist a fixed distribution of $X$ and a fixed concept $C$ such that the expected error is larger than a constant times $k/n$ for infinitely many $n$, where $k$ is the number of parameters. They obtained strong minimax lower bounds for the tail distribution of the probability of error, which extend the corresponding minimax lower bounds.

Our second form of lower bound, that is, Theorem 1.2, is comparable with the individual lower rate of convergence in [1]. In the latter, the individual lower rate of convergence for a class $D$ of distributions of $(X,Y)$ is defined by $a_n$ which satisfies

$$\inf_{\hat{g}_n} \sup_{P\in D} \limsup_{n\to\infty} a_n^{-1}\Big(L_P(\hat{g}_n) - \min_g L_P(g)\Big) > 0, \tag{1.13}$$

where $g$ is a classifier and $\hat{g}_n$ is an estimator. A class of distributions $D_\beta$ of $(X,Y)$ is given as the product of one uniform distribution and a cubic class of regression functions with parameter $\beta$. Under some assumptions, the individual lower rate of convergence for $D_\beta$ is obtained by $b_n n^{-2\beta/(2\beta+d)}$. The class $D_\beta$ is of course different from our class $\mathcal{P}$, but the order of inf, sup and lim sup in (1.13) is the same as ours in Theorem 1.2.

For more references on lower bounds, see Gill and Levit [8] and Tsybakov [18].

For a more general nonparametric setup than ours, Pfanzagl showed in [14] that *no* limit distribution can be attained with the rate $n^{1/3}$ *uniformly* on certain shrinking neighborhoods of the sample distribution $P$.

In a paper by Bühlmann and Yu in [4], the $n^{1/3}$-asymptotic appears in the context of bagging. These authors also use the results of Kim and Pollard [9]. Using decision trees, problems concerning higher dimensional $X$ are reduced to the analysis of a one-dimensional setup.

*Organization of this paper.* Let us explain how the rest of this paper is organized. In Section 2 we collect some fundamental entropy estimates. Section 3 shows *universal*, general counterparts of Theorems 1.1 and 1.2, without assuming the specific form of our model $(\Omega, \mathcal{A}, \mathcal{P})$. We expect these lemmas to be useful for other examples too. One key idea is the use of Baire's theorem to show that the set of $P$'s with estimation errors being asymptotically large on a given scale is of the second Baire category. In Section 4 Theorems 1.1 and 1.2 are proven. The proofs are based on a bound for relative entropies for slightly perturbed densities, described in Lemma 4.1 below. In Section 5 Theorem 1.3 is shown by constructing the estimators $\hat{a}_{n,L}$ in a two-step procedure. Section 6 contains the proofs of the corollaries. The key idea for the higher dimensional case is sketched in Section 7.

**2. Preliminaries.** In this section we review some standard estimates to compare probabilities with respect to different measures, based on bounds of the relative entropy. Alternatively (and more or less equivalently), one could use bounds for the Hellinger distance instead of the relative entropy, but we do not follow this alternative approach here. Here, we need not assume any specific form of the model $(\Omega,\mathcal{A},\mathcal{P})$; we take an arbitrary parameter $a:\mathcal{P}\to\mathbb{R}$, and $(\hat{a}_n:\Omega^n\to\mathbb{R})_{n\in\mathbb{N}}$ denotes any sequence of estimators.

Let $H(P,Q):=E_P[\log\frac{dP}{dQ}]$ denote the relative entropy for $P,Q\in\mathcal{P}$, whenever it is well defined.

LEMMA 2.1. *Let $P$ and $Q$ be probability measures with $H(P,Q)<\infty$. For every random variable $X$ with $0\le X\le 1$, one has*

$$E_Q[X]\ge e^{-2H(P,Q)-1}(E_P[X]-\tfrac{1}{2}). \tag{2.1}$$

PROOF. For $x\ge 0$, set $\psi(x):=x\log x-x+1$. Note that $\psi\ge 0$. We set $N:=e^{2H(P,Q)+1}\ge e$ and $A:=\{dP/dQ>N\}$. Using $\psi\ge 0$, $\psi(1)=0$ and the convexity of $\psi$, one sees that

$$\psi(x)\ge\frac{\psi(N)}{N}\mathbb{1}\{x>N\}x \tag{2.2}$$

for all $x\ge 0$, and thus,

$$H(P,Q)=E_Q\left[\psi\left(\frac{dP}{dQ}\right)\right]\ge\frac{\psi(N)}{N}E_Q\left[\mathbb{1}(A)\frac{dP}{dQ}\right]=\frac{\psi(N)}{N}P[A]. \tag{2.3}$$

We conclude, using $\psi(N)/N=\log(N/e)+1/N\ge 2H(P,Q)$,

$$\begin{aligned}E_Q[X]&\ge E_Q[X\mathbb{1}(A^c)]\ge\frac{1}{N}E_P[X\mathbb{1}(A^c)]\ge\frac{1}{N}(E_P[X]-P[A])\\&\ge\frac{1}{N}\left(E_P[X]-\frac{N}{\psi(N)}H(P,Q)\right)\ge\frac{1}{N}\left(E_P[X]-\frac{1}{2}\right),\end{aligned} \tag{2.4}$$

which is the claim (2.1). □

The 2 in the exponent of (2.1) could be replaced by any fixed number larger than 1, if one replaced the $\frac{1}{2}$ in (2.1) by a different constant. This would only change the constants in our main theorems.

LEMMA 2.2. *Let $\chi:\mathbb{R}\to[0,1]$ be a measurable function with $\chi(x)=1$ for $|x|\le 1$ and $\chi(x)=0$ for $|x|\ge 2$. Take $n\in\mathbb{N}$, $\beta_n>0$, $P_n,Q_n\in\mathcal{P}$, and $\delta\in(0,1/4]$. If*

$$nH(P_n,Q_n)\le\frac{1}{2}\log\frac{1}{11\delta},\qquad\beta_n|a(P_n)-a(Q_n)|>4 \tag{2.5}$$

*hold, then at least one of the two bounds*

$$\begin{aligned}&E_{P_n^n}[\chi(\beta_n(\hat{a}_n-a(P_n)))]<1-\delta\quad\textit{or}\\&E_{Q_n^n}[\chi(\beta_n(\hat{a}_n-a(Q_n)))]<1-\delta\end{aligned} \tag{2.6}$$

*is valid.*

PROOF. (Indirectly). Assume that both formulas in (2.6) fail to hold. Using Lemma 2.1, we get

$$\begin{aligned}
&Q_n^n[|\beta_n(\hat a_n - a(P_n))| \le 2] \\
&\quad \ge E_{Q_n^n}[\chi(\beta_n(\hat a_n - a(P_n)))] \\
&\quad \ge e^{-2nH(P_n,Q_n)-1}\left(E_{P_n^n}[\chi(\beta_n(\hat a_n - a(P_n)))] - \frac{1}{2}\right) \\
&\quad \ge e^{-2nH(P_n,Q_n)-1}\left(1 - \delta - \frac{1}{2}\right) \ge \frac{1}{4e} e^{-2nH(P_n,Q_n)} \ge \frac{11}{4e}\delta > \delta
\end{aligned} \tag{2.7}$$

by (2.5); recall that $\delta \le 1/4$. Furthermore, we have

$$Q_n^n[\beta_n|\hat a_n - a(Q_n)| \le 2] \ge E_{Q_n^n}[\chi(\beta_n(\hat a_n - a(Q_n)))] \ge 1 - \delta \tag{2.8}$$

by the opposite of the right-hand side of (2.6) and the choice of $\chi$. Recall that $\chi : \mathbb{R} \to [0,1]$. The bounds (2.8) and (2.7) imply that the events $\{\beta_n|\hat a_n - a(Q_n)| \le 2\}$ and $\{\beta_n|\hat a_n - a(P_n)| \le 2\}$ have a nonempty intersection. This implies the contradiction $\beta_n|a(P_n) - a(Q_n)| \le 4$. $\square$

**3. General lower bounds.** In this section we prepare the proofs of Theorems 1.1 and 1.2 by providing some general, abstract lower bounds for estimators. Since we expect these lemmas to be useful also in contexts other than the estimation of thresholds, we do not assume $\Omega$, $\mathcal{P}$ and $a : P \mapsto a(P)$ have the specific form described in Section 1. Rather, in this section $(\Omega, \mathcal{A})$ may be any measurable space, $\mathcal{P}$ may be any set of probability measures on $(\Omega, \mathcal{A})$, endowed with a topology, and $a : \mathcal{P} \to \mathbb{R}$ may be any parameter; only some very general assumptions for the topology on $\mathcal{P}$ and for the parameter $a$ are required in Lemma 3.2 below. The first lemma in this section is a general statement which plays an essential role in the proof of Theorem 1.1.

LEMMA 3.1. *Take a sequence $(\beta_n)_{n\in\mathbb{N}}$ of positive numbers, $\delta \in (0, 1/4]$, and a nonempty open set $\mathcal{U} \subseteq \mathcal{P}$. For all large $n \in \mathbb{N}$, assume that there are $P_n, Q_n \in \mathcal{U}$ such that*

$$nH(P_n, Q_n) \le \frac{1}{2}\log\frac{1}{11\delta} \quad \text{and} \quad \beta_n|a(P_n) - a(Q_n)| > 4 \tag{3.1}$$

*hold. Then for all sequences $(\hat a_n : \Omega^n \to \mathbb{R})_{n\in\mathbb{N}}$ of estimators, one has*

$$\liminf_{n\to\infty} \sup_{P\in\mathcal{U}} P^n[\beta_n|\hat a_n - a(P)| > 1] \ge \delta. \tag{3.2}$$

PROOF. Take $\chi := \mathbb{1}[-1,1]$. By (3.1) and Lemma 2.2, (2.6) holds for all large $n$. Thus,

$$\inf_{P\in\mathcal{U}} P^n[\beta_n|\hat a_n - a(P)| \le 1] < 1 - \delta \tag{3.3}$$

holds for all large $n$; hence, (3.2) is true. $\square$

The constant 4 on the right-hand side in (3.1) is not optimal. However, improving it does not improve our main theorems.

The next lemma provides an abstract key ingredient for the proof of Theorem 1.2.

LEMMA 3.2. *Let $a:\mathcal{P}\to\mathbb{R}$ be a parameter and $(\hat{a}_n:\Omega^n\to\mathbb{R})_{n\in\mathbb{N}}$ be a sequence of estimators. Let $\mathcal{P}$ be endowed with a Baire space metrizable topology. Assume that $a:\mathcal{P}\to\mathbb{R}$ is continuous. Furthermore, assume that for all $P\in\mathcal{P}$, the total variation distance from $P$, that is, $Q\mapsto\|P-Q\|_{\mathcal{A}}=\sup_{A\in\mathcal{A}}|P[A]-Q[A]|$, $Q\in\mathcal{P}$, is continuous too. Let $(\beta_n)_{n\in\mathbb{N}}$ be a sequence of positive numbers, and take $\delta\in(0,1/4]$.*

*Suppose that for all $P\in\mathcal{P}$, for all neighborhoods $\mathcal{N}$ of $P$, and for all $m\in\mathbb{N}$, there are $n\geq m$ and $Q_n\in\mathcal{N}$ such that*

$$\beta_n|a(P)-a(Q_n)|>4 \quad \text{and} \quad nH(P,Q_n)\leq\tfrac{1}{2}|\log(11\delta)| \tag{3.4}$$

*are valid. Then for all nonempty open sets $\mathcal{U}\subseteq\mathcal{P}$, one has*

$$\sup_{P\in\mathcal{U}}\limsup_{n\to\infty}P^n[\beta_n|\hat{a}_n-a(P)|>1]\geq\delta. \tag{3.5}$$

PROOF. Let $\chi:\mathbb{R}\to[0,1]$ be a continuous function with $\chi(x)=1$ for $|x|\leq 1$, and $\chi(x)=0$ for $|x|\geq 2$. For $m,n\in\mathbb{N}$, we set

$$\mathcal{P}_n:=\{P\in\mathcal{P}:E_{P^n}[\chi(\beta_n(\hat{a}_n-a(P)))]\geq 1-\delta\}, \qquad \mathcal{F}_m:=\bigcap_{n\geq m}\mathcal{P}_n. \tag{3.6}$$

We claim that the map

$$\mathcal{P}\to[0,1], \quad P\mapsto E_{P^n}[\chi(\beta_n(\hat{a}_n-a(P)))] \tag{3.7}$$

is continuous. To prove this claim, let $P\in\mathcal{P}$, and consider a sequence $(Q_k)_k$ in $\mathcal{P}$ converging to $P$. Then for all $\omega\in\Omega^n$, we have

$$\chi(\beta_n(\hat{a}_n(\omega)-a(Q_k)))\overset{k\to\infty}{\longrightarrow}\chi(\beta_n(\hat{a}_n(\omega)-a(P))) \tag{3.8}$$

by the continuity of $a$ and of $\chi$. Using Lebesgue's dominated convergence theorem, this implies

$$E_{P^n}[\chi(\beta_n(\hat{a}_n-a(Q_k)))]\overset{k\to\infty}{\longrightarrow}E_{P^n}[\chi(\beta_n(\hat{a}_n-a(P)))]; \tag{3.9}$$

recall that $\chi$ takes values in the unit interval. Furthermore,

$$\begin{aligned}&|E_{Q_k^n}[\chi(\beta_n(\hat{a}_n-a(Q_k)))]-E_{P^n}[\chi(\beta_n(\hat{a}_n-a(Q_k)))]|\\&\quad\leq n\|Q_k-P\|_{\mathcal{A}}\overset{k\to\infty}{\longrightarrow}0,\end{aligned} \tag{3.10}$$

since, by our hypothesis, the total variation distance from $P$ is continuous. Combining (3.9) and (3.10), we get

$$E_{Q_k^n}[\chi(\beta_n(\hat{a}_n - a(Q_k)))] \stackrel{k\to\infty}{\longrightarrow} E_{P^n}[\chi(\beta_n(\hat{a}_n - a(P)))], \tag{3.11}$$

which shows that $E_{P^n}[\chi(\beta_n(\hat{a}_n - a(P)))]$ depends continuously on $P$. Note that in the last step we used the fact that the chosen topology on $\mathcal{P}$ is metrizable (or, at least, that sequential continuity on $\mathcal{P}$ implies continuity).

The continuity of the map described in (3.7) implies that the sets $\mathcal{P}_n \subseteq \mathcal{P}$ are closed; thus, their intersections $\mathcal{F}_m$ are closed too.

Next, we show that the sets $\mathcal{F}_m \subseteq \mathcal{P}$, $m \in \mathbb{N}$, are nowhere dense. To check this, take $P \in \mathcal{F}_m$ and a neighborhood $\mathcal{N}$ of $P$ in $\mathcal{P}$. By the hypothesis of the lemma, there exist $n \geq m$ and $Q_n \in \mathcal{N}$ such that (3.4) holds. Then Lemma 2.2 implies

$$\begin{aligned} &E_{P^n}[\chi(\beta_n(\hat{a}_n - a(P)))] < 1-\delta \quad \text{or} \\ &E_{Q_n^n}[\chi(\beta_n(\hat{a}_n - a(Q_n)))] < 1-\delta, \end{aligned} \tag{3.12}$$

that is, $P \notin \mathcal{P}_n \supseteq \mathcal{F}_m$ or $Q_n \notin \mathcal{P}_n \supseteq \mathcal{F}_m$, and thus, $\mathcal{N} \not\subseteq \mathcal{F}_m$. This shows that indeed $\mathcal{F}_m$ is nowhere dense.

Let $\mathcal{U} \subseteq \mathcal{P}$ be a nonempty open set. Since $\mathcal{P}$ is endowed with a Baire space topology, we conclude that $\mathcal{U}$ is not contained in $\bigcup_{m\in\mathbb{N}} \mathcal{F}_m$; so we can take $P \in \mathcal{U} \setminus \bigcup_{m\in\mathbb{N}} \mathcal{F}_m$. For this $P$, we know that $P \notin \mathcal{P}_n$ for infinitely many $n \in \mathbb{N}$. Thus, we get

$$\liminf_{n\to\infty} E_{P^n}[\chi(\beta_n(\hat{a}_n - a(P)))] \leq 1-\delta. \tag{3.13}$$

Using

$$P^n[\beta_n|\hat{a}_n - a(P)| > 1] \geq 1 - E_{P^n}[\chi(\beta_n(\hat{a}_n - a(P)))], \tag{3.14}$$

this implies

$$\limsup_{n\to\infty} P^n[\beta_n|\hat{a}_n - a(P)| > 1] \geq \delta, \tag{3.15}$$

and thus, the claim (3.5) follows. □

**4. Lower bounds for errors in threshold estimation.** In this section we prove Theorems 1.1 and 1.2. Thus, let $\mathcal{P}$ denote again the concrete space of probability measures defined in Section 1, endowed with the metric $d$, defined in (1.3). Note that $(\tilde{\mathcal{P}}, d)$ is a *complete* metric space. We claim that $a: \mathcal{P} \to [0,1]$ is a *continuous* parameter. Indeed, this is a consequence of the implicit function theorem applied to the map $F: (0,1) \times C^1([0,1] \times \{\pm 1\}) \to \mathbb{R}$, $F(x,f) := f(x,1) - f(x,-1)$. Theorem 10.2.1 in [5] presents a version of the implicit function theorem applicable to our situation. The map $F$ is continuously differentiable with the derivative

$$\begin{aligned} DF(x,f)\colon&(\Delta x,\Delta f) \\ &\mapsto [D_1 f(x,1)-D_1 f(x,-1)]\Delta x+\Delta f(x,1)-\Delta f(x,-1); \end{aligned} \tag{4.1}$$

thus, the implicit function theorem is applicable for any point $(x,f)$ for which the transversality condition $D_1 f(x,1)\neq D_1 f(x,-1)$ holds. It yields the continuity of the function $a\colon\mathcal{P}\to(0,1)$, implicitly defined by the equation $F(a(P),f_P)=0$. Furthermore, $\mathcal{P}$ is an *open* subset of the space $\tilde{\mathcal{P}}$. In particular, by Baire's category theorem, $\mathcal{P}$ is a Baire space.

The following lemma contains the basic entropy estimate for perturbed densities. Here is the idea. A given probability density $f_P$ is slightly modified by a perturbation of order $O(\varepsilon)$ in a neighborhood of size $O(\varepsilon)$ of the transversal intersection point $a(P)$. Let $Q$ denote the probability measure corresponding to the modified density $f_Q$; then we show that the entropy $H(P,Q)$ has roughly the order $O(\varepsilon^3)$, but the parameter $a(Q)$ deviates from $a(P)$ on a scale of order $\varepsilon$. The cube of $\varepsilon$ arising in the entropy bound is the key to derive the cube root asymptotic lower bounds in this paper.

LEMMA 4.1. *Let $P\in\mathcal{P}$, and let $\mathcal{U}\subseteq\mathcal{P}$ be an open neighborhood of $P$. Then there is $c_1=c_1(P,\mathcal{U})>0$ such that for every $\delta\in(0,1/4]$ and for all large $n$ [say for $n\geq n_0(P,\mathcal{U},\delta)$], there is $Q_n\in\mathcal{U}$ such that*

$$nH(P,Q_n)\leq\tfrac{1}{2}|\log(11\delta)| \quad \textit{and} \quad \frac{n^{1/3}}{c_1|\log(11\delta)|^{1/3}}|a(P)-a(Q_n)|>4. \tag{4.2}$$

PROOF. Choose a ball $\mathcal{N}\subseteq\mathcal{U}$ (with respect to the metric $d$) centered at $P$. Let $r$ denote the radius of $\mathcal{N}$. Let $c_1>0$ be small enough (to be specified below). Take $\delta\in(0,1/4]$. We abbreviate

$$\beta_n := \frac{n^{1/3}}{c_1|\log(11\delta)|^{1/3}}>0. \tag{4.3}$$

Take a fixed, compactly supported $\phi\in C^1(\mathbb{R})$ with $\phi(0)=1$, $0\leq\phi\leq1$, and $\|\phi'\|_\infty\|f_P\|_\infty<r$. For $\varepsilon>0$, we set

$$\Xi_\varepsilon(x):=\varepsilon\phi\biggl(\frac{x-a(P)}{\varepsilon}\biggr). \tag{4.4}$$

Recall the definitions of $f_P^+$ and $f_P^-$ in (1.1). We set

$$\rho_P^+:=\frac{f_P^+}{f_P^++f_P^-},\qquad \rho_P^-:=\frac{f_P^-}{f_P^++f_P^-}; \tag{4.5}$$

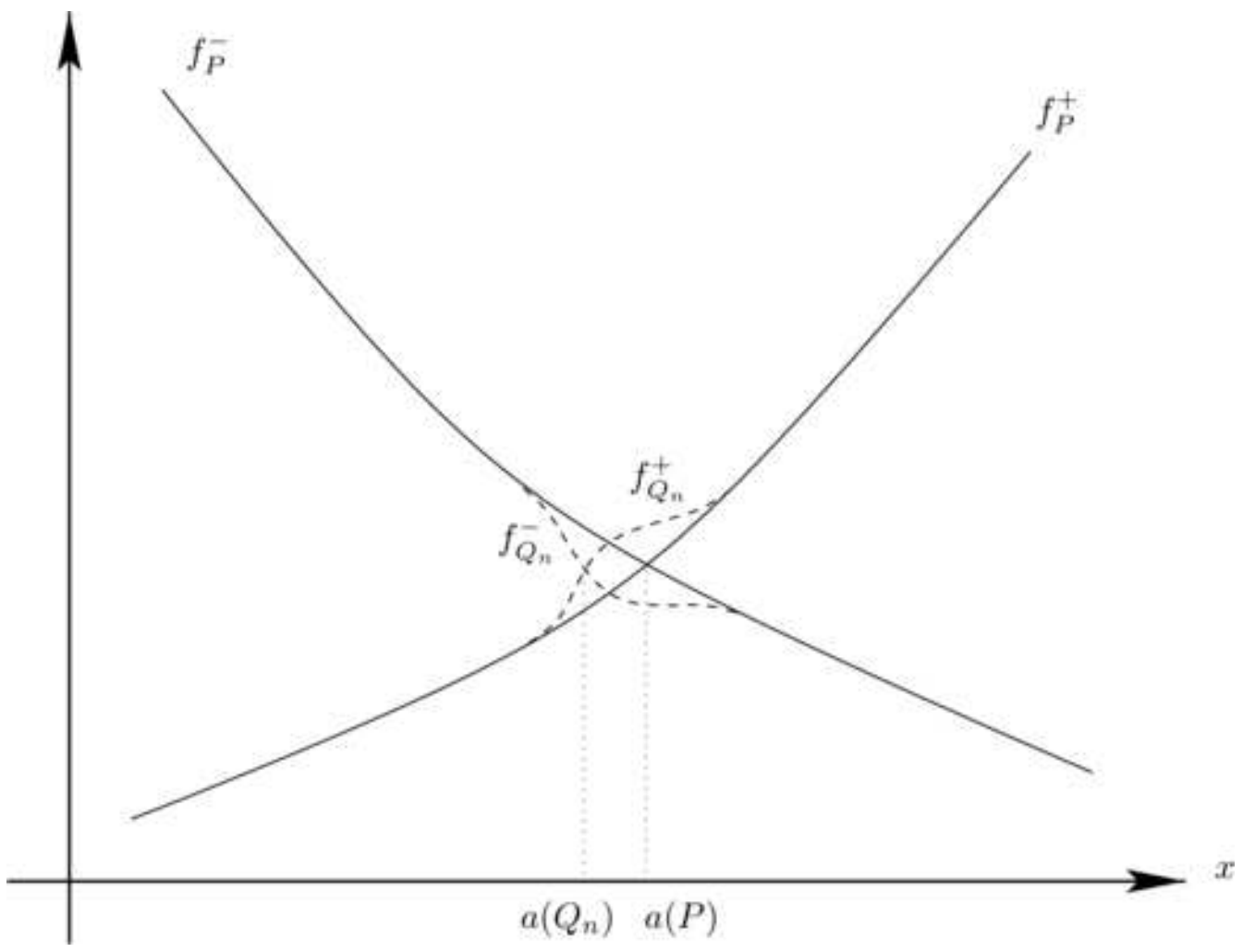

FIG. 1. *Perturbation of the two densities.*

at least in some compact neighborhood $V_P$ of $a(P)$, these functions are well defined with values in $[1/3, 2/3]$. For all small $\varepsilon > 0$, $\Xi_\varepsilon$ is supported in such a neighborhood $V_P$. We set

$$\varepsilon_n := c_4 |\log(11\delta)|^{1/3} n^{-1/3}, \tag{4.6}$$

where $c_4 = c_4(P, \phi) := (\|f_P\|_\infty \|\phi\|_2^2)^{-1/3}$. Let $Q_n$ be defined by its density $f_{Q_n}$, where

$$f_{Q_n}^\pm := (1 + \Xi_{\varepsilon_n} \rho_P^\mp) f_P^\pm, \tag{4.7}$$

$$f_{Q_n}(x, y) := f_{Q_n}^+(x) \mathbb{1}\{y = 1\} + f_{Q_n}^-(x) \mathbb{1}\{y = -1\}. \tag{4.8}$$

Figure 1 illustrates these definitions.

Here $\Xi_{\varepsilon_n} \rho_P^\pm$ is to be interpreted as 0 outside the support of $\Xi_{\varepsilon_n}$. Note that for all large $n$, $f_{Q_n}^\pm$ is well defined. As a consequence of the assumptions (1.4), one sees that $f_P^+(a(P)) = f_P^-(a(P)) > 0$. For large $n$, $\varepsilon_n$ is small; thus $f_{Q_n}^\pm$ is nonnegative and $f_{Q_n}$ is a probability density. Furthermore, using $\|\phi'\|_\infty \|f_P\|_\infty < r$ and $|\rho_P^\mp f_P^\pm| \le \|f_P\|_\infty$, one sees that $d(Q_n, P) < r$ and, thus, $Q_n \in \mathcal{N} \subseteq \mathcal{U}$ holds for all large $n$. We calculate, for $(x, y) \in \Omega$,

$$\frac{dQ_n}{dP}(x, y) = 1 + (\mathbb{1}\{y = 1\} \rho_P^-(x) - \mathbb{1}\{y = -1\} \rho_P^+(x)) \Xi_{\varepsilon_n}(x). \tag{4.9}$$

For $|t - 1| \le 1/2$, one has $-\log t + t - 1 \le |t - 1|^2$. So

$$H(P, Q_n) = E_P\left[-\log \frac{dQ_n}{dP} + \frac{dQ_n}{dP} - 1\right] \le E_P\left[\left(\frac{dQ_n}{dP} - 1\right)^2\right] \tag{4.10}$$

for $|\frac{dQ_n}{dP}-1|\leq 1/2$, which holds for all large $n$. Note that $|\Xi_{\varepsilon_n}|\leq\varepsilon_n$. For all large $n$, $\rho_P^\pm\in[1/3,2/3]$ holds on the support of $\Xi_{\varepsilon_n}$. Then one has

$$\tfrac{1}{3}\leq|\mathbb{1}\{y=1\}\rho_P^-(x)-\mathbb{1}\{y=-1\}\rho_P^+(x)|\leq\tfrac{2}{3} \tag{4.11}$$

on the support of $\Xi_{\varepsilon_n}$, and therefore,

$$\left(\frac{dQ_n}{dP}-1\right)^2\leq\left(\frac{2}{3}\right)^2\Xi_{\varepsilon_n}^2(X). \tag{4.12}$$

We get the following $O(\varepsilon_n^3)$ estimate for the relative entropy:

$$\begin{aligned}nH(P,Q_n)&\leq\tfrac{4}{9}nE_P[\Xi_{\varepsilon_n}^2(X)]\leq\tfrac{4}{9}n\|f_P\|_\infty\|\Xi_{\varepsilon_n}\|_2^2\\&\leq\tfrac{1}{2}n\|f_P\|_\infty\|\phi\|_2^2\varepsilon_n^3\leq\tfrac{1}{2}|\log(11\delta)|\end{aligned} \tag{4.13}$$

by the choice (4.6) of $\varepsilon_n$.

[As a side remark, note that the estimate (4.13) relies on our choice to take the Lebesgue measure $\lambda_{[0,1]}$ in the reference measure. In some cases with arbitrary reference measures it would break down.]

On the other hand, defining $\rho_{Q_n}^\pm$ in analogy to (4.5), from (4.7) and using $f_{Q_n}^+ + f_{Q_n}^- = f_P^+ + f_P^-$, it follows that

$$\rho_{Q_n}^\pm:=(1+\Xi_{\varepsilon_n}\rho_P^\mp)\rho_P^\pm. \tag{4.14}$$

Using $\rho_P^-(a(P))=\rho_P^+(a(P))=\rho_{Q_n}^+(a(Q_n))=1/2$ and (4.14), we get

$$\begin{aligned}\frac{\varepsilon_n}{4}&=\frac{1}{4}|\Xi_{\varepsilon_n}(a(P))|=|\Xi_{\varepsilon_n}(a(P))\rho_P^-(a(P))\rho_P^+(a(P))|\\&=|\rho_P^+(a(P))-\rho_{Q_n}^+(a(P))|=|\rho_{Q_n}^+(a(Q_n))-\rho_{Q_n}^+(a(P))|\\&\leq\|(\rho_{Q_n}^+)'\|_{\infty,V_P}|a(Q_n)-a(P)|.\end{aligned} \tag{4.15}$$

Taking the derivative of (4.14) and taking a supremum over $V_P$, we see that $\|(\rho_{Q_n}^+)'\|_{\infty,V_P}$ is bounded by a constant $c_5(P,\phi)>0$ for all large $n$; note that $\|\Xi_{\varepsilon_n}'\|_\infty=\|\phi'\|_\infty$ does not depend on $n$. We obtain

$$\begin{aligned}\beta_n|a(Q_n)-a(P)|&\geq\beta_n\frac{\varepsilon_n}{4c_5(P)}\\&=\frac{n^{1/3}}{c_1|\log(11\delta)|^{1/3}}\frac{c_4|\log(11\delta)|^{1/3}n^{-1/3}}{4c_5(P,\phi)}>4\end{aligned} \tag{4.16}$$

when we choose

$$0<c_1(P,\mathcal{U})<\frac{c_4(P,\phi)}{16c_5(P,\phi)}; \tag{4.17}$$

recall the choices (4.3) and (4.6) of $\beta_n$ and $\varepsilon_n$. The statements (4.16) and (4.13) together are just the claim (4.2). □

PROOF OF THEOREM 1.1. Take a fixed $P \in \mathcal{U}$, and take $c_1(\mathcal{U}) = c_1(P,\mathcal{U}) > 0$ from Lemma 4.1. Then Lemma 4.1 guarantees that the hypothesis (3.1) of Lemma 3.1 holds with $P_n = P$, where $\beta_n$ is again given by (4.3). Thus Lemma 3.1 yields the claim (1.5). $\square$

PROOF OF THEOREM 1.2. The class $\mathcal{P}$ with the metric $d$ is indeed a Baire space, and $a : \mathcal{P} \to (0,1)$ is continuous. Note that the total variation distance is continuous with respect to $d$.

We check that Lemma 3.2 is applicable with $\delta = 1/4$. Let $P \in \mathcal{P}$, and let $\mathcal{N}$ be a neighborhood of $P$ in $\mathcal{P}$. We apply Lemma 4.1 to obtain a sequence $(Q_n)_n$ in $\mathcal{P}$ such that (4.2) holds. Hence, we get, for all large $n$ [say, for $n \geq c_6(P,\mathcal{U})$],

$$\beta_n|a(P) - a(Q_n)| \geq \frac{n^{1/3}}{c_1(P,\mathcal{U})|\log(11\delta)|^{1/3}}|a(P) - a(Q_n)| > 4 \tag{4.18}$$

by $n^{-1/3}\beta_n \to \infty$. Together with the entropy bound in (4.2), this shows that Lemma 3.2 is indeed applicable, and it yields the claim (1.6). $\square$

Note that this proof would break down if we had taken $\beta_n$ on the *critical* scale $\beta_n = \text{const} \cdot n^{1/3}$. Indeed, the constant $c_1(P,\mathcal{U})$ depends on $\mathcal{U}$ (and it really diverges as $\mathcal{U}$ gets smaller), but $\beta_n$ must not depend on the choice of $\mathcal{U}$. This breakdown has to occur, since Theorem 1.3 shows that the claim (1.6) of Theorem 1.2 cannot hold any more on the critical scale $\beta_n = \text{const} \cdot n^{1/3}$.

**5. Optimal estimators for thresholds.** In this section we prove Theorem 1.3. The optimal estimators, whose errors asymptotically meet the border permitted by the lower bounds, are constructed by a two-step procedure. In the first step (Lemma 5.1 below) we use the empirical risk minimizer to obtain a *starting estimator*, which yields error terms roughly on the scale $n^{-1/3}$. The second step (Lemma 5.2 below) constructs a refined family of estimators, based on a starting approximation $a_0$. Let $(X_i, Y_i) : \Omega^n \to \Omega$ again denote the canonical projections.

LEMMA 5.1. *There is a sequence of estimators $\hat{a}_n : \Omega^n \to [0,1]$, $n \in \mathbb{N}$, such that, for all $P \in \mathcal{P}$, there is a neighborhood $\mathcal{N}(P)$ of $P$ in $\mathcal{P}$, such that*

$$\lim_{L\to\infty} \sup_{Q\in\mathcal{N}(P)} \lim_{n\to\infty} Q^n[|\hat{a}_n - a(Q)| \geq Ln^{-1/3}] = 0. \tag{5.1}$$

PROOF. Consider the empirical risk minimizers $\hat{a}_n$, $n \in \mathbb{N}$. We abbreviate $f_Q^\Sigma := f_Q^+ + f_Q^-$. Take $P \in \mathcal{P}$ and

$$\begin{aligned} \mathcal{N}(P) := \{Q \in \mathcal{P} : 2|(\rho_Q^+)'(a(Q))| > |(\rho_P^+)'(a(P))|, \\ 2f_Q^\Sigma(a(Q)) > f_P^\Sigma(a(P))\}, \end{aligned} \tag{5.2}$$

where $\rho_Q^+$ is again defined as in (4.5). By the transversality of the intersection point $a(Q)$ of $f_Q^+$ and $f_Q^-$, the maps $Q \mapsto (\rho_Q^+)'(a(Q))$ and $Q \mapsto f_Q^\Sigma(a(Q))$, $Q \in \mathcal{N}(P)$, are continuous. Furthermore, $f_P^\Sigma(a(P)) > 0$ and $|(\rho_P^+)'(a(P))| > 0$ hold. Using these facts, one sees that $\mathcal{N}(P)$ is a neighborhood of $P$.

Take $Q \in \mathcal{N}(P)$. By Theorem 2.2 in [13], we know that

$$n^{1/3}(\hat{a}_n - a(Q)) \xrightarrow{\mathcal{L}} [(\rho_Q^+)'(a(Q))\sqrt{f_Q^\Sigma(a(Q))}]^{-2/3} Z \tag{5.3}$$

as $n \to \infty$ for some continuous random variable $Z$ not depending on $Q$ (with respect to some probability measure $\mathbb{P}$). We set

$$\alpha := \inf_{Q \in \mathcal{N}(P)} [(\rho_Q^+)'(a(Q))\sqrt{f_Q^\Sigma(a(Q))}]^{2/3}. \tag{5.4}$$

Note that $\alpha > 0$ by the choice (5.2) of $\mathcal{N}(P)$. We now obtain the claim (5.1):

$$\begin{aligned} &\sup_{Q \in \mathcal{N}(P)} \lim_{n \to \infty} Q[|\hat{a}_n - a(Q)| \geq L n^{-1/3}] \\ &\quad = \sup_{Q \in \mathcal{N}(P)} \mathbb{P}[Z \geq [(\rho_Q^+)'(a(Q))\sqrt{f_Q^\Sigma(a(Q))}]^{2/3} L] \\ &\quad \leq \mathbb{P}[Z \geq \alpha L] \xrightarrow{L \to \infty} 0. \end{aligned} \tag{5.5}$$

□

The next lemma is used to construct the refined estimators in Theorem 1.3. Here is the idea. Given a starting approximation $a_0$ for $a(Q)$ with an error on the scale $n^{-1/3}$, consider all data points in a neighborhood of size $Ln^{-1/3}$, where $L$ is large, but fixed. Then construct a regression line through the data points in this neighborhood, and take the intersection of this regression line with the $x$-axis as the refined estimator.

LEMMA 5.2. *There is a family of estimators* $(\hat{a}_{n,L,a_0} : \Omega^n \to \mathbb{R})_{n \in \mathbb{N}, L > 0, a_0 \in (0,1)}$ *with the following property. For all* $P \in \mathcal{P}$, *there is a neighborhood* $\mathcal{N} \subseteq \mathcal{P}$ *of* $P$, *such that, for all* $T > 0$, *one has*

$$\lim_{L \to \infty} \sup_{Q \in \mathcal{N}} \limsup_{n \to \infty} \sup_{\substack{a_0 \in (0,1) \\ |a_0 - a(Q)| \leq L n^{-1/3}}} Q^n[n^{1/3}|\hat{a}_{n,L,a_0} - a(Q)| > T] = 0. \tag{5.6}$$

PROOF. Given $n \in \mathbb{N}$, $L > 0$, and $a_0 \in (0,1)$, we introduce the abbreviations $\tilde{X}_i := X_i - a_0$, $M := Ln^{-1/3}$ and $I_M := [-M, M]$, and we define the random set

$$J := \{j : 1 \leq j \leq n, \tilde{X}_j \in I_M\}. \tag{5.7}$$

We define the estimator $\hat{a}_{n,L,a_0}$ as follows: Consider the regression line $\ell$ (equation $y = \hat{b}_1 x + \hat{b}_2$) through the points $(\tilde{X}_j, Y_j)$, $j \in J$, provided it is well defined, that is, provided there are at least two different $\tilde{X}_{j_1} \neq \tilde{X}_{j_2}$, $j_1, j_2 \in J$. If $\hat{b}_1 \neq 0$, we set

$$\hat{a}_{n,L,a_0} := a_0 - \frac{\hat{b}_2}{\hat{b}_1}. \tag{5.8}$$

Geometrically this means that $\hat{a}_{n,L,a_0}$ is the intersection of the real axis with the regression line through the points $(X_j, Y_j)$, where only the points with $|X_j - a_0| \leq Ln^{-1/3}$ are taken, whenever this intersection is well defined. If the regression line $\ell$ is not well defined, or if $\hat{b}_1 = 0$, we set $\hat{a}_{n,L,a_0} = a_0$, just to have a definite value in this case too.

We abbreviate, for $Q \in \mathcal{P}$,

$$s_Q := f_Q^{\Sigma}(a(Q)), \qquad t_Q := (f_Q^+)'(a(Q)) - (f_Q^-)'(a(Q)), \tag{5.9}$$

where again $f_Q^{\Sigma} := f_Q^+ + f_Q^-$. Let $P \in \mathcal{P}$, and take the following neighborhood of $P$:

$$\mathcal{N} := \left\{ Q \in \mathcal{P} : \begin{array}{l} |t_Q| > |t_P|/2, s_Q < 2s_P, d(P,Q) < 1, \\ a(Q) > a(P)/2, 1 - a(Q) > (1 - a(P))/2 \end{array} \right\}. \tag{5.10}$$

Take $T > 0$ and $\delta > 0$. Let $L$ be large enough [more specifically, so large that

$$2L \geq T, \qquad L^3 \geq S^2 := \frac{5c_7}{\delta}, \qquad \frac{S}{\sqrt{L}} < \frac{T}{5} \tag{5.11}$$

hold with some positive constants $c_7 = c_7(\mathcal{N})$ and $c_8 = c_8(\mathcal{N})$, to be specified below]. We claim that, for all $Q \in \mathcal{N}$, one has

$$\limsup_{n \to \infty} \sup_{\substack{a_0 \in (0,1) \\ |a_0 - a(Q)| \leq Ln^{-1/3}}} Q^n[n^{1/3}|\hat{a}_{n,L,a_0} - a(Q)| > T] \leq \delta. \tag{5.12}$$

The claim (5.6) of the lemma then follows immediately from the statement (5.12).

Here is a sketch of the proof of (5.12). For a complete proof; see [12].

To prove (5.12), let $Q \in \mathcal{N}$. For every $n \in \mathbb{N}$ and for every $a_0 \in (0,1)$ with

$$|a_0 - a(Q)| \leq Ln^{-1/3}, \tag{5.13}$$

we are going to define an event $B = B(Q, n, a_0, S, L) \subseteq \Omega^n$ with the property

$$Q^n[B(Q, n, a_0, S, L)] \geq 1 - \delta. \tag{5.14}$$

[This is done in (5.22) below.] We then show that, for all large $n$ [say, for $n \geq n_0(Q, \mathcal{N}, S, L)$, but uniformly in the choice of $a_0$], one has

$$B(Q, n, a_0, S, L) \subseteq \{n^{1/3}|\hat{a}_{n,L,a_0} - a(Q)| \leq T\}. \tag{5.15}$$

Once we have proven (5.14) and (5.15), the claim (5.12) is an immediate consequence.

Take $Q \in \mathcal{N}$ and $a_0 \in (0,1)$ with the constraint (5.13). It is convenient to shift quantities by $a_0$. We set

$$\tilde{a} := \tilde{a}(Q) = a(Q) - a_0, \tag{5.16}$$

$$\tilde{f}_Q^{\pm}(x) := f_Q^{\pm}(x + a_0), \qquad \tilde{f}_Q^{\Sigma}(x) := f_Q^{\Sigma}(x + a_0). \tag{5.17}$$

In particular, note that $|\tilde{a}| \le M = Ln^{-1/3}$ and that for all large $n$ [uniformly in $a_0$ and $Q \in \mathcal{N}$ with the constraint (5.13)], $\tilde{f}_Q^{\pm}$ is defined at least on $I_M$.

The coefficients $\hat{b}_1$ and $\hat{b}_2$ of the regression line $\ell$ are determined by the linear system $A\hat{b} = c$, where

$$\begin{aligned} A &:= \begin{pmatrix} \sum_{j\in J} \tilde{X}_j^{\,2} & \sum_{j\in J} \tilde{X}_j \\ \sum_{j\in J} \tilde{X}_j & |J| \end{pmatrix}, \\ \hat{b} &:= \begin{pmatrix} \hat{b}_1 \\ \hat{b}_2 \end{pmatrix}, \\ c &:= \begin{pmatrix} \sum_{j\in J} \tilde{X}_j Y_j \\ \sum_{j\in J} Y_j \end{pmatrix}. \end{aligned} \tag{5.18}$$

We introduce the (normalized) difference of the coefficient matrix from its expected value,

$$\begin{pmatrix} M^{5/2}\Delta_{11} & M^{3/2}\Delta_{12} \\ M^{3/2}\Delta_{21} & M^{1/2}\Delta_{22} \end{pmatrix} = \frac{1}{\sqrt{n}}(A - E_{Q^n}(A)), \tag{5.19}$$

$$\begin{pmatrix} M^{3/2}\Gamma_1 \\ M^{1/2}\Gamma_2 \end{pmatrix} = \frac{1}{\sqrt{n}}(c - E_{Q^n}[c]). \tag{5.20}$$

Our reason to normalize the $\Delta_{ij}$ and $\Gamma_i$ in this way is the following bound on the variances (to be read element-wise):

$$\begin{aligned} \operatorname{Var}_{Q^n}\left[\begin{pmatrix} \Delta_{11} & \Delta_{12} \\ \Delta_{21} & \Delta_{22} \end{pmatrix}\right] &\le E_{Q^n}\left[\mathbb{1}\{\tilde{X}_1 \in I_M\} \begin{pmatrix} \tilde{X}_1^4/M^5 & \tilde{X}_1^2/M^3 \\ \tilde{X}_1^2/M^3 & 1/M \end{pmatrix}\right] \\ &\le \begin{pmatrix} c_7 & c_7 \\ c_7 & c_7 \end{pmatrix}, \\ \operatorname{Var}_{Q^n}\left[\begin{pmatrix} \Gamma_1 \\ \Gamma_2 \end{pmatrix}\right] \le E_{Q^n}&\left[\mathbb{1}\{\tilde{X}_1 \in I_M\} \begin{pmatrix} \tilde{X}_1^2 Y_1^2/M^3 \\ Y_1^2/M \end{pmatrix}\right] \le \begin{pmatrix} c_7 \\ c_7 \end{pmatrix}, \end{aligned} \tag{5.21}$$

with some constant $c_7 = c_7(\mathcal{N}) > 0$. We have used the fact that the density $f_Q$ of $(X_i, Y_i)$ is uniformly bounded for $Q \in \mathcal{N}$.

Here is the definition of the event $B$:

$$B(Q,n,a_0,S,L) := \{|\Delta_{i,j}| \le S, |\Gamma_i| \le S \ (i,j \in \{1,2\})\}. \tag{5.22}$$

Chebyshev's inequality, (5.21) and the choice (5.11) of $S$ imply the claim (5.14),

$$Q^n[B(Q,n,a_0,S,L)^c] \le \frac{5c_7}{S^2} = \delta. \tag{5.23}$$

The factor 5 arises since there are five random variables involved (recall $\Delta_{12} = \Delta_{21}$).

In the rest of this proof we verify the claim (5.15) for all large $n$ (uniformly in $a_0$). So assume that the event $B(Q,n,a_0,S,L)$ holds.

The system $A\hat{b} = c$ is equivalent to

$$\begin{aligned} &\begin{pmatrix} \int_{I_M} x^2 \tilde{f}_Q^\Sigma(x)\,dx + \frac{M^{5/2}\Delta_{11}}{\sqrt{n}} & \int_{I_M} x \tilde{f}_Q^\Sigma(x)\,dx + \frac{M^{3/2}\Delta_{12}}{\sqrt{n}} \\ \int_{I_M} x \tilde{f}_Q^\Sigma(x)\,dx + \frac{M^{3/2}\Delta_{21}}{\sqrt{n}} & \int_{I_M} \tilde{f}_Q^\Sigma(x)\,dx + \frac{M^{1/2}\Delta_{22}}{\sqrt{n}} \end{pmatrix} \hat{b} \\ &= n \begin{pmatrix} \int_{I_M} x(\tilde{f}_Q^+(x) - \tilde{f}_Q^-(x))\,dx + \frac{M^{3/2}\Gamma_1}{\sqrt{n}} \\ \int_{I_M} (\tilde{f}_Q^+(x) - \tilde{f}_Q^-(x))\,dx + \frac{M^{1/2}\Gamma_2}{\sqrt{n}} \end{pmatrix}. \end{aligned} \tag{5.24}$$

Let us introduce some notation used in the Taylor approximations below. The variables $\xi_j = \xi_j(x,\tilde{a},Q,a_0)$ denote some values between $x$ and $\tilde{a}$. The variables $\varepsilon_j$ denote error terms which are bounded by $\mathcal{N}$-dependent constants $|\varepsilon_j| \le \mathrm{const}(\mathcal{N})$, and $\delta_j$ denote error terms which are bounded by $|\delta_j| \le \mathrm{const} \cdot \sigma(Q,2M)$, where

$$\sigma(Q,r) := \max_{y=\pm 1} \sup_{|x_1-x_2|\le r} |\partial_1 f_Q(x_1,y) - \partial_1 f_Q(x_2,y)| \overset{r\to 0}{\longrightarrow} 0 \tag{5.25}$$

denotes the modulus of continuity of $\partial_1 f_Q$; recall that $f_Q \in C^1(\Omega)$. Recall that $|\tilde{a}| \le M$. Let us approximate the integrals in (5.24) by Taylor's formula,

$$\begin{aligned} \int_{I_M} x^2 \tilde{f}_Q^\Sigma(x)\,dx &= \int_{-M}^{M} x^2 (s_Q + (x-\tilde{a})(\tilde{f}_Q^\Sigma)'(\xi_1))\,dx \\ &= \tfrac{2}{3} M^3 s_Q + \varepsilon_1 M^4, \end{aligned} \tag{5.26}$$

$$\int_{I_M} x \tilde{f}_Q^\Sigma(x)\,dx = \int_{-M}^{M} x(s_Q + (x-\tilde{a})(\tilde{f}_Q^\Sigma)'(\xi_2))\,dx = \varepsilon_2 M^3, \tag{5.27}$$

$$\int_{I_M} \tilde{f}_Q^\Sigma(x)\,dx = \int_{-M}^{M} (s_Q + (x-\tilde{a})(\tilde{f}_Q^\Sigma)'(\xi_3))\,dx$$

$$= 2Ms_Q + \varepsilon_3 M^2, \tag{5.28}$$

$$\begin{aligned}\int_{I_M} x(\tilde f_Q^+(x) - \tilde f_Q^-(x))\,dx &= \int_{-M}^{M} x(x-\tilde a)((\tilde f_Q^+)'(\xi_4) - (\tilde f_Q^-)'(\xi_4))\,dx \\ &= \tfrac{2}{3}M^3 t_Q + \delta_1 M^3 = \varepsilon_4 M^3,\end{aligned} \tag{5.29}$$

$$\begin{aligned}\int_{I_M} (\tilde f_Q^+(x) - \tilde f_Q^-(x))\,dx &= \int_{-M}^{M} (x-\tilde a)((\tilde f_Q^+)'(\xi_5) - (\tilde f_Q^-)'(\xi_5))\,dx \\ &= -2\tilde a M t_Q + \delta_2 M^2 = \varepsilon_5 M^2.\end{aligned} \tag{5.30}$$

We rewrite the system (5.24) as

$$\begin{aligned}&\begin{pmatrix} \frac{2}{3}M^3 s_Q + \varepsilon_1 M^4 + \frac{M^{5/2}\Delta_{11}}{\sqrt{n}} & \varepsilon_2 M^3 + \frac{M^{3/2}\Delta_{12}}{\sqrt{n}} \\ \varepsilon_2 M^3 + \frac{M^{3/2}\Delta_{21}}{\sqrt{n}} & 2Ms_Q + \varepsilon_3 M^2 + \frac{M^{1/2}\Delta_{22}}{\sqrt{n}} \end{pmatrix} \hat b \\ &= \begin{pmatrix} \frac{2}{3}M^3 t_Q + \delta_1 M^3 + \frac{M^{3/2}\Gamma_1}{\sqrt{n}} \\ -2\tilde a M t_Q + \delta_2 M^2 + \frac{M^{1/2}\Gamma_2}{\sqrt{n}} \end{pmatrix} \\ &= \begin{pmatrix} \varepsilon_4 M^3 + \frac{M^{3/2}\Gamma_1}{\sqrt{n}} \\ \varepsilon_5 M^2 + \frac{M^{1/2}\Gamma_2}{\sqrt{n}} \end{pmatrix};\end{aligned} \tag{5.31}$$

both forms of the right-hand side are useful below. Dividing the first row in (5.31) by $(2/3)M^3 s_Q$ and the second row by $2Ms_Q$, we get the normalized system

$$\begin{aligned}&\begin{pmatrix} 1 + M\left[\frac{3}{2}\frac{\varepsilon_1}{s_Q} + \frac{3}{2s_Q}\frac{\Delta_{11}}{L^{3/2}}\right] & \frac{3}{2}\frac{\varepsilon_2}{s_Q} + \frac{3}{2s_Q}\frac{\Delta_{12}}{L^{3/2}} \\ M^2\left[\varepsilon_6 + \frac{1}{2s_Q}\frac{\Delta_{21}}{L^{3/2}}\right] & 1 + M\left[\varepsilon_7 M + \frac{1}{2s_Q}\frac{\Delta_{22}}{L^{3/2}}\right] \end{pmatrix} \hat b \\ &= \frac{t_Q}{s_Q}\begin{pmatrix} 1 + \frac{3}{2}\frac{\delta_1}{t_Q} + \frac{3}{2t_Q}\frac{\Gamma_1}{L^{3/2}} \\ -\tilde a + M\left(\frac{\delta_2}{2t_Q} + \frac{1}{2t_Q}\frac{\Gamma_2}{L^{3/2}}\right) \end{pmatrix} \\ &= \frac{t_Q}{s_Q}\begin{pmatrix} \varepsilon_8 + \varepsilon_9\frac{\Gamma_1}{L^{3/2}} \\ M\left[\varepsilon_{10} + \varepsilon_{11}\frac{\Gamma_2}{L^{3/2}}\right] \end{pmatrix}.\end{aligned} \tag{5.32}$$

Heuristically, one should view (5.32) as a perturbation of the system

$$\begin{pmatrix} 1 & 0 \\ 0 & 1 \end{pmatrix} \begin{pmatrix} b_1 \\ b_2 \end{pmatrix} = \frac{t_Q}{s_Q} \begin{pmatrix} 1 \\ -\tilde{a} \end{pmatrix}, \tag{5.33}$$

for which one knows $-b_2/b_1 = \tilde{a}$.

By (5.11) and the definition (5.22) of the event $B$, we know, for $i, j \in \{1,2\}$, $\frac{|\Delta_{i,j}|}{L^{3/2}}, \frac{|\Gamma_i|}{L^{3/2}} \le \frac{S}{L^{3/2}} \le 1$ and $|\frac{3}{2t_Q} \frac{\Gamma_i}{L^{3/2}}| \le c_8 \frac{S}{L^{3/2}} \le \frac{T}{5L}$, when we choose the constant in (5.11) to be $c_8 = c_8(\mathcal{N}) = \sup_{Q \in \mathcal{N}} |3/(2t_Q)| \le 3|t_P|^{-1} < \infty$; recall the definition (5.10) of $\mathcal{N}$, and recall that we assume the event $B$ holds.

For $M \le 1$, we rewrite (5.32) in the form

$$\begin{aligned} \begin{pmatrix} 1 + M\varepsilon_{12} & \varepsilon_{13} \\ M^2 \varepsilon_{14} & 1 + M\varepsilon_{15} \end{pmatrix} \hat{b} &= \frac{t_Q}{s_Q} \begin{pmatrix} 1 + \varepsilon_{16}\delta_1 + \varepsilon_{17}\dfrac{T}{L} \\ -\tilde{a} + M\left[\varepsilon_{18}\delta_2 + \varepsilon_{19}\dfrac{T}{L}\right] \end{pmatrix} \\ &= \frac{t_Q}{s_Q} \begin{pmatrix} \varepsilon_{20} \\ M\varepsilon_{21} \end{pmatrix}, \end{aligned} \tag{5.34}$$

where $|\varepsilon_{19}| < 1/5$ and $|\varepsilon_{17}| \le 1/5$.

Let us consider the asymptotics of $\hat{b}_2/\hat{b}_1$ as $n \to \infty$, that is, as $M = Ln^{-1/3} \to 0$. For all large $n$ (uniformly in $a_0$), the system (5.34) is nonsingular; recall that the error terms $\varepsilon_j$ are bounded uniformly in $a_0$. We get, for all large $n$,

$$\frac{\hat{b}_2}{\hat{b}_1} = \frac{\begin{vmatrix} 1 + M\varepsilon_{12} & \varepsilon_{20} \\ M^2\varepsilon_{14} & -\tilde{a} + M[\varepsilon_{18}\delta_2 + \varepsilon_{19}T/L] \end{vmatrix}}{\begin{vmatrix} 1 + \varepsilon_{16}\delta_1 + \varepsilon_{17}T/L & \varepsilon_{13} \\ M\varepsilon_{21} & 1 + M\varepsilon_{15} \end{vmatrix}}. \tag{5.35}$$

Let $d_2$ and $d_1$ denote the determinants in the numerator and denominator of the right-hand side in (5.35), respectively. Recall that the error terms $\delta_2$, $\delta_1$ converge to 0 as $n \to \infty$, uniformly in $a_0$; see (5.25). Using $|\varepsilon_{19}|, |\varepsilon_{17}| \le 1/5 < 1/4$, we get, for all large $n$ (uniformly in $a_0$),

$$|-\tilde{a} - d_2| \le \frac{MT}{4L}, \qquad |1 - d_1| \le \frac{T}{4L}, \qquad \left|1 - \frac{1}{d_1}\right| \le \frac{T}{2L}. \tag{5.36}$$

We have used $T \le 2L$ in the last step.

Using the definition (5.8) of $\hat{a}_{n,L,a_0}$ and (5.35), we conclude, still for large $n$ (uniformly in $a_0$),

$$\begin{aligned} |\hat{a}_{n,L,a_0} - a(Q)| &= \left| -\tilde{a} - \frac{\hat{b}_2}{\hat{b}_1} \right| = \left| -\tilde{a} - \frac{d_2}{d_1} \right| \\ &\le |\tilde{a}| \frac{T}{2L} + \frac{MT}{4L} + \frac{MT}{4L}\frac{T}{2L} \le \frac{MT}{L} = Tn^{-1/3}; \end{aligned} \tag{5.37}$$

recall that $|\tilde{a}| \le M$ and $T \le 2L$. Thus, the claim (5.15) holds for all large $n$, uniformly in $a_0$. $\square$

We now use Lemma 5.1 to obtain the starting approximation $a_0$ required by Lemma 5.2.

Proof of Theorem 1.3. Let us abbreviate $Z_j := (X_j, Y_j)$. Let $n \in \mathbb{N}$, $n \ge 2$, and $L > 0$. We construct the estimator $\hat{a}_{n,L}$ by a two-step procedure. We split the sample $Z_1, \dots, Z_n$ into two halves $Z_1, \dots, Z_m$ and $Z_{m+1}, \dots, Z_{2m}$, where we abbreviate $m = m(n) := \lfloor n/2 \rfloor$. (For odd $n$, we drop one data point at the end.) We then use the first half of the data to get a *rough* estimate $\hat{a}_m = \hat{a}_m(Z_1, \dots, Z_m)$ by Lemma 5.1. Using this as a *starting estimate*, we refine it by Lemma 5.2, applied to the second half of the data

$$\hat{a}_{n,L} := \hat{a}_{m,L,\hat{a}_m}(Z_{m+1}, \dots, Z_{2m}). \tag{5.38}$$

It is important to split the data into two disjoint pieces, since then $\hat{a}_m$ is independent of $(Z_{m+1}, \dots, Z_{2m})$; thus, we have good control of the distribution of

$$\hat{a}_{m,L,\hat{a}_m}(Z_{m+1}, \dots, Z_{2m})$$

conditioned on $\hat{a}_m(Z_1, \dots, Z_m)$. By a slight abuse of notation, we abbreviate

$$\hat{a}_{m,L,a_0} := \hat{a}_{m,L,a_0}(Z_{m+1}, \dots, Z_{2m}).$$

Let $P \in \mathcal{P}$, and let $\mathcal{N}$ be the intersection of the two neighborhoods of $P$ that were constructed in Lemmas 5.1 and 5.2. Let $T > 0$ and $\delta > 0$. By the same two lemmas, we know, for all large $L$,

$$\sup_{Q \in \mathcal{N}(P)} \lim_{n \to \infty} Q^n[A(Q,n,L)^c] < \frac{\delta}{2} \tag{5.39}$$

with the events

$$A(Q,n,L) := \{|\hat{a}_m - a(Q)| \le L m^{-1/3}\} \tag{5.40}$$

and

$$\sup_{Q \in \mathcal{N}} \limsup_{n \to \infty} \sup_{\substack{a_0 \in (0,1) \\ |a_0 - a(Q)| \le L m^{-1/3}}} Q^n\left[m^{1/3} |\hat{a}_{m,L,a_0} - a(Q)| > \frac{T}{2}\right] < \frac{\delta}{2}. \tag{5.41}$$

Note that, for all large $n$ [say, for $n \ge n_0(Q,L)$], we have $\hat{a}_m \in (0,1)$ on the event $A(Q,n,L)$. Let $Q \in \mathcal{N}$, and let $n \in \mathbb{N}$ be large enough. Using

$n^{1/3} \le 2m^{1/3}$ in the first step, we get,

$$
\begin{aligned}
&Q^n[n^{1/3}|\hat{a}_{n,L} - a(Q)| > T] \\
&\quad \le Q^n\left[m^{1/3}|\hat{a}_{n,L} - a(Q)| > \frac{T}{2}\right] \\
&\quad \le E_{Q^n}\left[Q^n\left[m^{1/3}|\hat{a}_{n,L} - a(Q)| > \frac{T}{2}\Big|\hat{a}_m\right]\mathbb{1}(A(Q,n,L))\right] \qquad (5.42)\\
&\qquad + Q^n[A(Q,n,L)^c] \\
&\quad \le E_{Q^n}\left[Q^n\left[m^{1/3}|\hat{a}_{m,L,\hat{a}_m} - a(Q)| > \frac{T}{2}\Big|\hat{a}_m\right]\mathbb{1}(A(Q,n,L))\right] + \frac{\delta}{2}.
\end{aligned}
$$

Using the independence structure and (5.41), we know for almost all $a_0$ with $|a_0 - a(Q)| \le Lm^{-1/3}$ (*almost all* with respect to the law of $\hat{a}_m$)

$$
\begin{aligned}
&Q^n\left[m^{1/3}|\hat{a}_{m,L,\hat{a}_m} - a(Q)| > \frac{T}{2}\Big|\hat{a}_m = a_0\right] \\
&\qquad = Q^n\left[m^{1/3}|\hat{a}_{m,L,a_0} - a(Q)| > \frac{T}{2}\right] \le \frac{\delta}{2}.
\end{aligned} \qquad (5.43)
$$

Combining (5.42) and (5.43), we conclude $Q^n[n^{1/3}|\hat{a}_{n,L} - a(Q)| > T] \le \delta$. This finishes the proof of the claim (1.7). □

Let us finally describe a simple counterexample, showing that the $\limsup_{n\to\infty}$ in the claim (1.6) of Theorem 1.2 cannot be replaced by $\liminf_{n\to\infty}$.

Let us take $\hat{a}_n$ to be *constant* estimators in the following way. For all $k \in \mathbb{N}_0$ and $n \in [2^k, 2^{k+1}[ \cap \mathbb{N}$, set $\hat{a}_n := (n - 2^k)2^{-k}$. Then, whatever the value $a(P) \in [0,1]$ is, we have the following: For all $k \in \mathbb{N}_0$, there is $n \in [2^k, 2^{k+1}[ \cap \mathbb{N}$ with $|\hat{a}_n - a(P)| \le 2^{-k} \le 2/n$. But then

$$\liminf_{n\to\infty} P^n[\beta_n|\hat{a}_n - a(P)| > 1] = 0 \qquad (5.44)$$

whenever $\beta_n = o(n)$ as $n \to \infty$.

Intuitively speaking, the counterexample uses the following idea. If you have many nonrunning clocks, one for every minute of the day, all showing different times, then one of them will show the correct time, up to 1 minute.

**6. Asymptotic bounds for the classification error.** In this section we present the proofs of Corollaries 1.4, 1.5 and 1.6. The proofs depend on the following lemma, which is based on a Taylor expansion.

LEMMA 6.1. *For all $Q \in \mathcal{P}$, there is a neighborhood $\mathcal{U}$ of $Q$ in $\mathcal{P}$, and there are positive constants $c_9 = c_9(\mathcal{U})$, $c_3 = c_3(\mathcal{U})$ and $c_{10} = c_{10}(\mathcal{U})$, such that, for all $P \in \mathcal{U}$ and for all $\alpha \in \mathbb{R}$*

$$\min\{c_9, c_3(a(P) - \alpha)^2\} \le L_P(\alpha) - L_P(a(P)) \le c_{10}(a(P) - \alpha)^2. \qquad (6.1)$$

PROOF. For $P\in\mathcal{P}$, we set $m_P := f_P^+ - f_P^-$. Note that $m_P(a(P)) = 0$. Furthermore, for all $\alpha\in(0,1)$,

$$L_P(\alpha) - L_P(a(P)) = \int_{a(P)}^{\alpha} m_P(x)\,dx = \int_{a(P)}^{\alpha} m'_P(x)(\alpha - x)\,dx. \tag{6.2}$$

Since $m'_P$ (in the $\|\cdot\|_\infty$-norm) and $a(P)$ depend continuously on $P$, the fact $m'_Q(a(Q)) > 0$ implies the following for some neighborhood $\mathcal{U}$ of $Q$ and some $\varepsilon = \varepsilon(\mathcal{U}) > 0$. For all $P\in\mathcal{U}$, one has $[a(P)-\varepsilon, a(P)+\varepsilon]\subset(0,1)$,

$$\begin{aligned} c_3 &:= \tfrac{1}{2}\inf_{P\in\mathcal{U}}\inf_{x\in[a(P)-\varepsilon,a(P)+\varepsilon]} m'_P(x) > 0 \quad\text{and}\\ c_{10} &:= \tfrac{1}{2}\sup_{P\in\mathcal{U}}\|m'_P\|_\infty < \infty. \end{aligned} \tag{6.3}$$

Thus, we get the upper bound in (6.1). Moreover, for all $\alpha\in(0,1)$ and $P\in\mathcal{U}$ with $|\alpha - a(P)|\le\varepsilon$, we have $c_3(a(P)-\alpha)^2 \le L_P(\alpha) - L_P(a(P))$. Since $\operatorname{sign}(m_P(x)) = \operatorname{sign}(x - a(P))$ holds for all $x\in[0,1]$, we have $L_P(\alpha)\ge L_P(a(P)-\varepsilon)$ for $\alpha < a(P)-\varepsilon$, and $L_P(\alpha)\ge L_P(a(P)+\varepsilon)$ for $\alpha > a(P)+\varepsilon$. Hence, we get the lower bound in (6.1) with $c_9 := c_3\varepsilon^2$. $\square$

PROOF OF COROLLARIES 1.4 AND 1.5. Consider an arbitrary sequence $\gamma_n$ with $\gamma_n \stackrel{n\to\infty}{\longrightarrow}\infty$. Then, it follows from the lower bound in (6.1) that, for $c_3T^2 < c_9\gamma_n^2$, that is, for all large $n$, we have

$$\{\gamma_n|\hat{a}_n - a(P)| > T\} \subseteq \{\gamma_n^2(L_P(\hat{a}_n) - L_P(a(P))) > c_3T^2\} \tag{6.4}$$

for all $P\in\mathcal{P}$. Take any $U$ and set $T := (U/c_3)^{1/2}$. In view of (6.4),

$$\begin{aligned} &\liminf_{n\to\infty}\sup_{P\in\mathcal{U}} P^n[\gamma_n^2(L_P(\hat{a}_n) - L_P(a(P))) > U]\\ &\quad\ge \liminf_{n\to\infty}\sup_{P\in\mathcal{U}} P^n[\gamma_n|\hat{a}_n - a(P)| > T]. \end{aligned}$$

Now, Corollary 1.4 follows from Theorem 1.1, taking $U := S$, $c_2 := c_1^2c_3$ and $\gamma_n := n^{1/3}$.

Similarly, Corollary 1.5 is a consequence of Theorem 1.2, by taking $U := c_3$ and $\gamma_n := \beta_n$. $\square$

PROOF OF COROLLARY 1.6. From the upper bound in (6.1), it follows that

$$\{n^{2/3}(L_P(\hat{a}_n) - L_P(a(P))) > T\} \subseteq \{n^{1/3}|\hat{a}_n - a(P)| > \sqrt{T/c_{10}}\}. \tag{6.5}$$

Now, we have

$$\begin{aligned} &\lim_{L\to\infty}\sup_{Q\in\mathcal{N}}\limsup_{n\to\infty} Q^n[n^{2/3}(L_Q(\hat{a}_n) - L_Q(a(Q))) > T]\\ &\quad\le \lim_{L\to\infty}\sup_{Q\in\mathcal{N}}\limsup_{n\to\infty} Q^n[n^{1/3}|\hat{a}_n - a(P)| > \sqrt{T/c_{10}}] = 0, \end{aligned} \tag{6.6}$$

where we have used the claim (1.7) of Theorem 1.3 in the last step. $\square$

**7. General theory for higher dimensions.** In this section we explain briefly how one can generalize our theory to higher dimensions. The full discussion of the generalization is beyond the scope of this paper. Consider a statistical model $\mathcal{P}$, that is, a class of probability measures over a measurable space $(\Omega, \mathcal{A})$. As an example, one could think of the law of a sample drawn according to two unknown smooth densities over a higher dimensional space, intersecting each other transversally in a hypersurface. Consider another measurable space $(\mathcal{H}, \mathcal{B}(\mathcal{H}))$ and a loss function

$$L : \mathcal{P} \times \mathcal{H} \to \mathbb{R}, \qquad (P, h) \mapsto L_P(h). \tag{7.1}$$

Assume that, for each $P \in \mathcal{P}$, there exists a minimizer $h_P \in \mathcal{H}$, that is,

$$L_P(h_P) \leq L_P(h) \qquad \text{for all } h \in \mathcal{H}. \tag{7.2}$$

Set $\Delta_P(h) := L_P(h) - L_P(h_P)$. For $\gamma > 0$, set $\Delta_P^\gamma(h) := \min\{\gamma, \Delta_P(h)\}$. The following lemma generalizes Lemma 2.2.

Lemma 7.1. *Let $P, Q \in \mathcal{P}$ and $\gamma > 0$ such that $\Delta_P(h) + \Delta_Q(h) \geq \gamma$, $\forall h \in \mathcal{H}$. Then for any $\delta \in (0, 1/2)$ and for any estimator $\hat{h}_n : \Omega^n \to \mathcal{H}$, at least one of the following two statements holds:*

$$E_{P^n}(\Delta_P^\gamma(\hat{h}_n)) \geq \delta\gamma \tag{7.3}$$

*or*

$$E_{Q^n}(\Delta_Q^\gamma(\hat{h}_n)) \geq (\tfrac{1}{2} - \delta)\gamma \exp(-2nH(P, Q) - 1). \tag{7.4}$$

Proof. Note that $\Delta_P^\gamma$ and $\Delta_Q^\gamma$ take values in $[0, \gamma]$. It is easily seen that $\Delta_P^\gamma(h) + \Delta_Q^\gamma(h) \geq \gamma$, for all $h \in \mathcal{H}$. Hence, $E_{P^n}(\Delta_P^\gamma(\hat{h}_n) + \Delta_Q^\gamma(\hat{h}_n)) \geq \gamma$, and for any $\delta \in (0, 1/2)$, we have $E_{P^n}(\Delta_P^\gamma(\hat{h}_n)) \geq \delta\gamma$ or $E_{P^n}(\Delta_Q^\gamma(\hat{h}_n)) \geq (1 - \delta)\gamma$. By Lemma 2.1, we know that

$$\begin{aligned} \frac{1}{\gamma} E_{Q^n}(\Delta_Q^\gamma(\hat{h}_n)) &\geq \left(\frac{1}{\gamma} E_{P^n}(\Delta_Q^\gamma(\hat{h}_n)) - \frac{1}{2}\right) \exp(-2nH(P, Q) - 1) \\ &= \frac{1}{\gamma}\left(E_{P^n}(\Delta_Q^\gamma(\hat{h}_n)) - \frac{\gamma}{2}\right) \exp(-2nH(P, Q) - 1). \end{aligned} \tag{7.5}$$

So, if $E_{P^n}(\Delta_Q^\gamma(\hat{h}_n)) \geq (1 - \delta)\gamma$, then

$$\begin{aligned} E_{Q^n}(\Delta_Q^\gamma(\hat{h}_n)) &\geq \left((1 - \delta)\gamma - \frac{\gamma}{2}\right) \exp(-2nH(P, Q) - 1) \\ &= (1/2 - \delta)\gamma \exp(-2nH(P, Q) - 1). \end{aligned} \tag{7.6}$$

□

In a higher dimensional setup, beyond the scope of this paper, Lemma 7.1 can be used as a replacement for Lemma 2.2. The assumption $\beta_n |a(P) -$ $a(Q_n)| > 4$ then becomes $\inf_h(\Delta_P^{\gamma_n}(h) + \Delta_{Q_n}^{\gamma_n}(h)) \geq \gamma_n$, where $\gamma_n$ depends on the problem. The lower bounds are obtained for the probability of $\Delta_P^{\gamma_n}(\hat{h}_n)$ being large, where $\hat{h}_n$ is an estimator.

**Acknowledgment.** We are grateful to D. Donoho, S. A. van de Geer and W. R. van Zwet for their support and for interesting discussions.

## REFERENCES

[1] Antos, A. (1999). Lower bounds on the rate of convergence of nonparametric pattern recognition. *Computational Learning Theory. Lecture Notes in Comput. Sci.* **1572** 241–252. Springer, Berlin. MR1724993

[2] Antos, A., Györfi, L. and Kohler, M. (2000). Lower bounds on the rate of convergence of nonparametric regression estimates. *J. Statist. Plann. Inference* **83** 91–100. MR1741445

[3] Antos, A. and Lugosi, G. (1998). Strong minimax lower bounds for learning. *Machine Learning* **30** 31–56.

[4] Bühlmann, P. and Yu, B. (2002). Analyzing bagging. *Ann. Statist.* **30** 927–961. MR1926165

[5] Dieudonné, J. (1969). *Foundations of Modern Analysis.* Academic Press, New York. MR0349288

[6] Donoho, D. L. and Liu, R. C. (1991). Geometrizing rates of convergence. II, III. *Ann. Statist.* **19** 633–667, 668–701. MR1105839

[7] Freund, Y. and Schapire, R. E. (1998). Large margin classification using the perceptron algorithm. In *Proc. Eleventh Annual Conference on Computational Learning Theory* 209–217. ACM, New York. MR1811585

[8] Gill, R. D. and Levit, B. Y. (1995). Applications of the Van Trees inequality: A Bayesian Cramér–Rao bound. *Bernoulli* **1** 59–79. MR1354456

[9] Kim, J. and Pollard, D. (1990). Cube root asymptotics. *Ann. Statist.* **18** 191–219. MR1041391

[10] Mammen, E. and Tsybakov, A. B. (1999). Smooth discrimination analysis. *Ann. Statist.* **27** 1808–1829. MR1765618

[11] Mason, L., Bartlett, P. L. and Golea, M. (2002). Generalization error of combined classifiers. *J. Comput. System Sci.* **65** 415–438. MR1939778

[12] Mohammadi, L. (2004). Estimation of thresholds in classification. Ph.D. dissertation, Leiden Univ.

[13] Mohammadi, L. and van de Geer, S. A. (2003). On threshold-based classification rules. In *Mathematical Statistics and Applications*: *Festschrift for Constance van Eeden* (M. Moore, S. Froda and C. Léger, eds.) 261–280. IMS, Beachwood, OH. MR2138297

[14] Pfanzagl, J. (2000). On local uniformity for estimators and confidence limits. *J. Statist. Plann. Inference* **84** 27–53. MR1747496

[15] Putter, H. and van Zwet, W. R. (1996). Resampling: Consistency of substitution estimators. *Ann. Statist.* **24** 2297–2318. MR1425955

[16] Schapire, R. E., Freund, Y., Bartlett, P. and Lee, W. S. (1998). Boosting the margin: A new explanation for the effectiveness of voting methods. *Ann. Statist.* **26** 1651–1686. MR1673273

[17] Stone, C. J. (1982). Optimal global rates of convergence for nonparametric regression. *Ann. Statist.* **10** 1040–1053. MR0673642

[18] Tsybakov, A. B. (2004). *Introduction à l'estimation non-paramétrique*. Springer, Berlin. MR2013911

[19] van de Geer, S. A. (2000). *Applications of Empirical Process Theory*. Cambridge Univ. Press. MR1739079

[20] van der Vaart, A. W. and Wellner, J. A. (1996). *Weak Convergence and Empirical Processes. With Applications to Statistics.* Springer, New York. MR1385671

Mathematical Institute
University of Munich
Theresienstrasse 39
D-80333 Munich
Germany
E-mail: merkl@mathematik.uni-muenchen.de

Medische Statistiek
LUMC, S5-P
Postbus 9600
NL-2300 RC Leiden
The Netherlands
E-mail: L.Mohammadi_Khankahdani@lumc.nl